\documentclass[runningheads]{cl2emult}

\usepackage{makeidx}  
\usepackage{graphicx} 
\usepackage{subeqnar} 
\usepackage{multicol} 
\usepackage{cropmark} 
\usepackage{math}     
\makeindex            

\usepackage{amsmath}
\usepackage{amssymb}
\usepackage{amsfonts}

\usepackage{xy}

\input xy
\xyoption{all}

\def\bC{{\mathbb C}}

\def\bQ{{\mathbb Q}}
\def\bZ{{\mathbb Z}}

\def\bP{{\mathbb P}}
\def\L{{\mathcal L}}

\def\bF{{\mathbb  F}}

\def\bG{{\bar G}}
\def\bg{{\bar g}}
\def\g{{g}}

\def\a{{\alpha}}

\def\sem{{\rtimes}}
\def\iso{{\, \cong\, }}
\def\e{{\xi}}
\def\r{{\bar r}}
\def\s{{\bar s}}
\def\z{{\bar z}}

\def\P{\mathcal P}
\def\<{\langle}
\def\>{\rangle}

\usepackage{amsmath}   
\begin{document}

\title*{Elliptic subfields and automorphisms of genus 2 function fields}

\author{T. Shaska
\and H. V\"olklein}


\institute{Department of Mathematics, University of Florida,
Gainesville, FL 32611.   }

\maketitle

\def\M{\mathcal M}

\begin{abstract}
We study genus 2 function fields with  elliptic subfields of degree
2. The locus $\L_2$ of these fields is a 2-dimensional subvariety of
the moduli space $\mathcal M_2$ of genus 2 fields.
An equation for $\L_2$ is already in the work of Clebsch and Bolza. We
use a birational parameterization of $\L_2$ by affine 2-space  
 to study the relation between the j-invariants of the degree 2 elliptic
subfields. This extends work of  Geyer, Gaudry, Stichtenoth and others.
We find a 1-dimensional family of genus 2 curves having exactly two 
isomorphic elliptic subfields of degree 2; this family is parameterized by the 
 j-invariant of these subfields.
\end{abstract}

{\it
\centerline{This paper is dedicated to Professor Shreeram Abhyankar}
\centerline{ on the occasion of his 70th birthday} 
}

\section{Introduction}
Sections 2 and 4 of this note are
 concerned with degree 2 elliptic subfields $E$ of a genus
2 function field $K$ (All function fields are over an algebraically
closed field $k$ of char. $\neq 2$). 
 Jacobi \cite{J} already noted that in this case $K$ has generators $X$
and $Y$ with
\begin{equation}\label{eq1}
Y^2=X^6-s_1X^4+s_2X^2-s_3
\end{equation}
This generalized an example of Legendre.  In the newer literature,
Cassels \cite{Cassels}  chapter 14 deals with arithmetic aspects of this.
Gaudry/Schost \cite{Gaudry} show that a  genus 2 field $K$ in $char > 5$
has at most two elliptic subfields of degree 2, up to isomorphism, and compute the
j-invariants of these elliptic subfields in terms of Igusa invariants
of $K$.

\par On the other hand, there is a group theoretic aspect.  
Degree 2 elliptic subfields of $K$ correspond to {\bf elliptic
involutions} in  the automorphism group of $K$  i.e. involutions different from the
hyperelliptic involution $e_0$. 
Thus our topic is intimately related with the structure of
$G:=Aut(K/k)$, and its quotient $\bG$ by $<e_0>$.
Geyer \cite{Geyer} 
classifies the possibilities for $\bG$, gives a brief discussion of $G$  
 and also notes some consequences for
isogenies between elliptic subfields. His exposition is very brief
because the main focus of his paper is on a different theme. We study
the structure of $G$ in section 3. We give a simple classification, based
on group-theoretic properties of central extensions of $\bG$,
and relate it to our $(u,v)$-parameterization of $\L_2$
(see below). It follows that the number of $G$-classes of
degree 2 elliptic subfields of $K$ is 0, 1 or 2; and this number is 1
if and only if $K$ has equation \ $Y^2=X(X^4-1)$.

 Brandt/Stichtenoth \cite{BS} more generally discuss automorphisms
of hyperelliptic curves (in characteristic  0), whereas Brandt \cite{Br}
(unpublished thesis) has a very comprehensive classification of
automorphism groups of hyperelliptic curves in any characteristic
and more generally, cyclic extensions of
genus 0 fields.

\par The purpose of this note is to combine these two aspects, the
geometric and the group theoretic one. E.g., Gaudry/Schost use only
the reduced automorphism group, using $G$ itself 
would  simplify their paper.
They  exclude characteristics 3 and 5 where  other
types of automorphism groups appear. 

\par In section 2 and 4 we study the locus $\L_2$ of
genus 2  fields with  elliptic subfields of degree
2. Geyer \cite{Geyer} states that $\L_2$ is a rational surface whose singular
locus is the curve corresponding to reduced automorphism group $V_4$
(see our section 3, case III). We give an explicit birational parametrization of $\L_2$
by parameters $u,v$; they are  obtained by setting $s_3=1$ in (1)
and symmetrizing $s_1,s_2$ by an action of $S_3$. More precisely,
those $u,v$ parametrize genus 2  fields together with an elliptic
involution of the reduced automorphism group (Thm \ref{sect1_thm}).
We express the j-invariants of degree 2 elliptic subfields in terms of $u,v$.
The particular case that these j-invariants are all equal (for a fixed
genus 2 field) yields a birational embedding of the moduli space
$\mathcal M_1$ of genus 1 curves into  $\mathcal M_2$. 

In section 4 we use the coordinates on $\M_2$ and $\L_2$ provided by
invariant theory. Expressing these coordinates in terms of our
$(u,v)$-parameters makes the parametrization of $\L_2$ explicit.
From this we confirm the explicit equation 
 found by Gaudry/Schost \cite{Gaudry} that is satisfied by all points
of $\L_2$; and we see directly that $\L_2$ is the full zero set
of this equation.

\par More generally, there is literature on degree $n$
elliptic subfields, e.g., Frey \cite{Frey},  and Frey and Kani
\cite{FK}, and Lange \cite{Lange}.   The first author's PhD thesis \cite{Sh} deals
with the case $n=3$. We further intend to study the cases
$n=5$ and 7. 

\par In the last section, we study the action of Aut$(K)$ on 
elliptic subfields $F$ of odd degree $n\ge 7$. 
The hyperelliptic involution fixes these subfields, hence they are permuted by $\bG$. 
It is asy to see that stabilizer $\bG_F$ in $\bG$ of F has order $\leq 3$. We study
 those cases where $\bG_F \neq 1$, assuming $char (k) =0$.
 This allows us to use Riemann's Existence Theorem to 
parametrize the extensions $K/F$ of degree $n$  with non-trivial automorphisms by
certain triples of permutations in $S_n$. To count the number of these
triples of permutations is a difficult problem for general $n$. We use
a computer search to construct all such triples for $n \le 21$.

\medskip\noindent
{\bf Notation:} All function fields in this paper are over $k$,
where  $k$ is an algebraically closed field of characteristic $\neq 2$.
Further, $V_4$ denotes the Klein 4-group and $D_{2n}$ (resp., $\bZ_n$)
the dihedral group of order $2n$ (resp., cyclic group of order $n$).

\section{Genus 2 Curves  with Elliptic Involutions}

\def\z{e} \def\e{\epsilon}
\def\A{\mbox{Aut}} \def\AA{\overline{\A}}

Let $K$ be a genus 2
field. Then $K$ has exactly one genus 0 subfield of degree 2, call it
$k(X)$. It is the fixed field of the 
{\bf hyperelliptic involution} $\z_0$ in $\A(K)$. Thus $\z_0$ is
central in $\A(K)$. Here and in the
following, $\A(K)$ denotes the group $\A(K/k)$, more precisely.
It induces a subgroup of $\A(k(X))$ which is naturally isomorphic
to $\AA(K):= \A(K)/<\z_0>$. The latter is called the
{\bf reduced automorphism group} of $K$.

\begin{definition} An {\bf elliptic involution} of $G = \A(K)$
 is  an  involution different from $\z_0$. Thus  the
elliptic involutions of $G$ are in  1-1 correspondence with the
elliptic subfields of $K$ of degree 2. An involution of 
$\bG=\AA(K)$ is called elliptic if it is the image of an
elliptic involution of $G$.
\end{definition}
 
If $\z_1$ is an elliptic involution in $G$
then $\z_2:=\z_0\, \z_1$  is another one.
So the elliptic involutions come naturally in
(unordered) pairs $\z_1$, $\z_2$.   
These pairs correspond bijectively to the elliptic involutions 
of $\bG$.
The latter also correspond to pairs $E_1, E_2$ of elliptic
subfields of $K$ of degree 2 with $E_1\cap k(X)=E_2\cap k(X)$.
 
\begin{definition}
We will consider pairs $(K,\e)$  with $K$ a genus 2 field and $\e$
an elliptic involution in $\bG$. Two such pairs $(K,\e)$ and  $(K',\e')$ are
called isomorphic if there is a $k$-isomorphism  $\a: K \to K'$ with
$\e' = \a \e \a^{-1}$.
\end{definition}

\par Let  $\e$ be an elliptic involution in $\bG$. We can choose the generator
$X$ of Fix$(\z_0)$      such that $\e(X)=-X$.  Then $K=k(X,Y)$ where 
$X, Y$ satisfy \eqref{eq1} 
with $s_1, s_2, s_3 \in k$, $s_3\ne0$ 
(follows from \eqref{2.0} and Remark \ref{rem0} in section 2). 
Further $E_1=k(X^2,Y)$
and $E_2=k(X^2,YX)$ are the two elliptic subfields corresponding to $\e$. 
Let $j_1$ and $j_2$ be their j-invariants.

\par Preserving the condition $\e(X)=-X$ we  can further modify $X$
such that  $s_3=1$. Then 
\begin{equation}\label{eq3}
Y^2=X^6-s_1X^4+s_2X^2-1
\end{equation}
where the polynomial on the right has non-zero discriminant. 

\def\zz{\zeta}

These conditions  determine $X$ up to coordinate
change by the group $\< \tau_1, \tau_2\>$ where  
$\tau_1:
X\to \zz_6X$, $\tau_2: X\to \frac 1 X$, and  $\zz_6$ is a primitive
6-th root of unity in $k$. (Thus $\zz_6=-1$ if char$(k)=3$).
 Here $\tau_1$ maps $(s_1,s_2)$ to
$(\zz_6^4s_1,\zz_6^2s_2)$, and $\tau_2$ switches $s_1,s_2$.
Invariants of this action are:

\begin{equation}\label{eq2}
\begin{split} 
u: & =s_1 s_2 \\
v: & =s_1^3+s_2^3
\end{split}
\end{equation}

In these parameters, the discriminant of the sextic
polynomial on the right hand side of \eqref{eq3}
equals $64 \Delta^2$, where  
$$\Delta\ =\ \Delta(u,v)\ =\ u^2-4v+18u-27\neq 0$$
Further, the j-invariants $j_1$ and $j_2$ are given by:

\begin{equation}\label{j_eq}
j_1+j_2 =256\frac {(v^2-2u^3+54u^2-9uv-27v)} {\Delta}
\end{equation}
$$j_1\, j_2 = 65536\frac {(u^2+9u-3v)} {\Delta^2}$$

The map $(s_1, s_2) \mapsto (u,v)$ is a branched Galois covering with
group $S_3$  of
the set $\{ (u,v)\in k^2 : \Delta(u,v)\neq 0\}$ by the corresponding
open subset of $s_1, s_2$-space  if char$(k)\ne3$. In any case, it is true
that if $s_1, s_2$ and $s_1', s_2'$ have the same $u,v$-invariants then
they are conjugate under $\< \tau_1, \tau_2\>$.

\begin{lemma}\label{lemma1}  For $(s_1,s_2)\in k^2$ with $\Delta\neq 0$, equation
\eqref{eq3} defines a  genus 2 field  $K_{s_1, s_2}=k(X,Y)$. Its
reduced automorphism group contains the elliptic involution
$\e_{s_1, s_2}: X \mapsto -X$.
Two such  pairs $(K_{s_1, s_2}, \e_{s_1, s_2})$ and  $(K_{s_1', s_2'},
\e_{s_1', s_2'})$ are isomorphic if and only if $u=u'$ and $v=v'$ 
(where $u,v$ and $u',v'$ are associated with $s_1, s_2$ and
$s_1', s_2'$, respectively, by \eqref{eq2}).
\end{lemma}

\proof
An isomorphism $\a$ between these two pairs yields $K=k(X,Y)=k(X', Y')$
with $ k(X)=k(X')$ such that $X,Y$ satisfy \eqref{eq3} and  $X',Y'$
satisfy the corresponding equation with $s_1, s_2$ replaced by $s_1',
s_2'$.
Further, $\e_{s_1, s_2}(X')=-X'$.
Thus $X'$ is conjugate to $X$ under 
 $\< \tau_1, \tau_2\>$ by the above remarks.
This proves the condition is necessary. It is clearly sufficient.

\begin{theorem}\label{sect1_thm}
i) The $(u,v)\in k^2$ with $\Delta\neq 0$ bijectively parameterize the
isomorphism classes of pairs $(K,\e)$ where $K$ is a genus 2 field and
$\e$ an elliptic involution of $\AA(K)$.
This parameterization is defined in Lemma
\ref{lemma1}. The j-invariants of the two elliptic subfields of $K$
associated with $\e$ are given by \eqref{j_eq}.

ii) The $(u,v)$ satisfying additionally

\begin{equation}\label{V_4}
(v^2-4u^3)(4v-u^2+110u-1125)\neq 0
\end{equation}

 bijectively parameterize the
isomorphism classes of genus 2 fields with $\A(K)\iso V_4$;
equivalently, genus 2 fields having exactly 2 elliptic subfields of
degree 2. Their j-invariants $j_1, j_2$ are given in terms of $u$ and $v$ by \eqref{j_eq}.
\end{theorem}

\proof
 i) follows from the Lemma. \\
iii) Condition \eqref{V_4} is equivalent to
$\A (K)$ being a Klein 4-group, and  to the other stated condition,
by 2.3, Case IV.
The theorem follows.

\begin{remark} (Isomorphic elliptic subfields) 
 For each $j\in k, j\neq 0, 1728, -32678$ there is
a unique genus 2 field $K$ with $\A(K)\cong V_4$ such that the two
elliptic subfields of $K$ of degree 2 have the same given j-invariant.
This generalizes as follows: For each $j\in k, j\neq 0$,
there is a   pair $(K,\e)$
as in the Theorem,  unique up to isomorphism, such that the two
associated elliptic subfields of $K$ have the same given
j-invariant and the corresponding $u,v$ satisfy $v=9(u-3)$. 
Mapping $j\in k\setminus \{0\}$ to the associated  $K$
gives an isomorphic embedding 
of $\mathcal M_1\setminus \{j=0\}$ into $\mathcal M_2$.  Here
$\mathcal M_g$ denotes the moduli space of genus $g$ curves (over
$k$). 
\end{remark}

\proof From \eqref{j_eq} we get that the discriminant of $(x-j_1)(x-j_2)$ is 
$$2^{16}\, (4u^3-v^2)(v-9u+27)^2 \Delta^2 $$
Thus the condition $j_1=j_2$ is equivalent to either
$v=9(u-3)$ or $v^2=4u^3$. The latter condition is equivalent to $\A
(K) \geq D_8$ by Lemma \ref{invol}(b) below.
 Under the condition $v=9(u-3)$ we get 
$$u= 9 - \frac {j} {256}, \quad v= 9(6- \frac {j} {256})$$
where $j:=j_1=j_2$.
There is only one point on the curve $v=9(u-3)$ with $\Delta (u,v)=0$,
namely $u=9$, $v=54$; it corresponds  to  $j=0$. Further, for
$j=1728$ (resp., $j=-32678$) we have $\A (K) \iso D_8$, (resp., $D_{12}$).
For all the
other values of $j$, we have $\A (K) \iso V_4$.
This proves the first claim  by part i). 
The rest is proved in section 3 using Igusa coordinates on
$\mathcal M_2$. 

\begin{remark} (2- and 3-isogenous elliptic subfields) 
The  modular 3-polynomial 
\begin{small}
\begin{equation}
\begin{split}
\Phi_3 &
=x^4-x^3y^3+y^4+2232xy(x+y)-1069956xy(x+y)+36864000(x^3+y^3)\\
& +2587918086x^2y^2+8900222976000xy(x+y)+452984832000000(x^2+y^2)\\
& -770845966336000000xy+1855425871872000000000(x+y)
\end{split}
\end{equation}
\end{small}
is symmetric in $j_1$ and $j_2$ hence becomes a polynomial in $u$ and $v$ via
\eqref{j_eq}. This polynomial factors as follows;
\begin{equation}
(4v-u^2+110u-1125)\cdot g_1(u,v)\cdot g_2(u,v)=0
\label{iso_3}
\end{equation}
where $g_1$ and $g_2$ are
\begin{small}
\begin{equation}
\begin{split}
g_1 & =-27008u^6+256u^7-2432u^5v+v^4+7296u^3v^2-6692v^3u-1755067500u\\
 & +2419308v^3-34553439u^4+127753092vu^2+16274844vu^3-1720730u^2v^2\\
 & -1941120u^5+381631500v+1018668150u^2-116158860u^3+52621974v^2\\
 & +387712u^4v -483963660vu-33416676v^2u+922640625 \\
\end{split}
\end{equation}
\end{small}
\begin{small}
\begin{equation}
\begin{split}
g_2 & =291350448u^6-v^4u^2-998848u^6v-3456u^7v+4749840u^4v^2+17032u^5v^2\\
 &  +4v^5+80368u^8+256u^9+6848224u^7-10535040v^3u^2-35872v^3u^3+26478v^4u\\
 &  -77908736u^5v+9516699v^4+307234984u^3v^2-419583744v^3u-826436736v^3\\
 &  +27502903296u^4+28808773632vu^2-23429955456vu^3+5455334016u^2v^2\\
 &  -41278242816v+82556485632u^2-108737593344u^3-12123095040v^2\\
 & +41278242816vu+3503554560v^2u+5341019904u^5-2454612480u^4v\\
\end{split}
\end{equation}
\end{small}
Vanishing of the first factor is equivalent to $D_{12} \leq G$, see part II 
of the next section. (Here again $G=\A(K)$). 
 If $G=D_{12}  $ then $K$ has two classes of
elliptic involutions $e$,  where  $e$ and $e_0 e$ are non-conjugate;
thus $K$ has two $G$-classes of elliptic subfields of degree 2,
and subfields from different classes are 3-isogenous. This was noted in
\cite{Gaudry} (for $p\neq 5$). 
There are exactly two fields $K$ such that 
$D_{12}$ is properly contained in $G$, see part I of the next section. In
these  cases, $e$ and $e_0 e$ are conjugate (and the corresponding
elliptic curves are 3-isogenous to themselves). 
 In the case III of the next section, $G$ has two classes of elliptic
involutions $e$; now $e$ and $e_0 e$ are conjugate, hence $j_1=j_2$ in formula
\eqref{j_eq}. Degree 2 elliptic subfields from different $G$-classes are
now 2-isogenous, see \cite{Geyer}.

\end{remark}

\section{Automorphism Groups of Genus 2 Fields}
\def\P{P} \def\a{p}

\subsection{Preliminaries} \label{3.1}

Let  $K$ be a genus 2 field, $G$ its automorphism group and
$\z_0\in G$ the hyperelliptic involution. Then 
$<\z_0>=Gal(K/k(X))$, where $k(X)$ is the unique genus $0$ subfield
of degree 2 of $K$. The reduced automorphism group 
$\bG=G/<\z_0>$ embeds into
$Aut(k(X)/k)\cong \mbox{PGL}_2(k)$. 

The  extension $K/k(X)$  is ramified at exactly  six places
$X=\a_1, \dots , \a_6$ of $k(X)$, where
$\a_1, \dots , \a_6$ are six distinct points in $\bP^1:= \bP^1_k$. Let
 $\P:=\{ \a_1, \dots , \a_6 \}$. The corresponding places of $K$ are called 
the {\bf Weierstrass points} of $K$. The group $G$ permutes the 6
Weierstrass points, and $\bG$ permutes accordingly $\a_1, \dots , \a_6$
in its action on $\bP^1$ as subgroup of $\mbox{PGL}_2(k)$. 
This yields an embedding $\bG \hookrightarrow S_6$. We have $K=k(X,Y)$, where

\begin{equation}\label{2.0}
Y^2=\prod_{\overset {\a\in \P}{ \a \neq \infty} }  (X-\a)
\end{equation}

\noindent
Because $K$ is the unique degree 2 extension of $k(X)$ ramified exactly at
$\a_1$, $\dots$ , $\a_6$, each automorphism of $k(X)$ permuting
these 6 places extends to  an automorphism of $K$.
Thus, $\bG$ is the stabilizer in $Aut(k(X)/k)\cong \mbox{PGL}_2(k)$ of the 6-set 
$\P$. 

Let $\Gamma:= \mbox{PGL}_2(k)$. If $l$ is
prime to $char(k)$ then  each  element of
order $l$ of $\Gamma$ is  conjugate to
  $\pmatrix {\e_l} & 0 \\ 0 & 1\endpmatrix$, where $\e_l$
is a  primitive $l$-th root of unity.  Each such element    has 2
fixed points on $\bP^1$ and
other orbits of length $l$.  If $l=char(k)$ then $\Gamma$
has exactly one class of elements of order $l$, represented by
$\pmatrix 1 & 1\\ 0 & 1\endpmatrix$. Each such element has exactly
one fixed point on $\bP^1$. 

\begin{lemma}\label{lem0}
Let $\g \in G$ and $\bg$ its image in $\bG$.

a)  Suppose $\bg$ is an involution. Then $\g$ has order 2 if and only
if it fixes no  Weierstrass points.

b) If $\bg$ has  order 4, then $\g$ has  order 8.
\end{lemma}
\proof  
a) Suppose $\bg $ is an involution. We may assume $\bg(X)=-X$.
\par Assume first that $\bg$ fixes no points in $\P$. Then
$\P=\{a,-a,b,-b,c,-c\}$ for certain $a,b,c \in k$. Thus
$$Y^2=(X^2-a^2)(X^2-b^2)(X^2-c^2)$$
and so $\g(Y)^2=Y^2$. 
Hence  $\g(Y)=\pm \, Y$, and $\g$ has order 2.
\par  Now suppose  $\bg$  fixes 2 points of $\P$. 
Then $\P=\{0,\infty, a,-a,b,-b\}$, hence
$$Y^2=X(X^2-a^2)(X^2-b^2)$$
So $\g(Y)^2=- Y^2      $ and $\g(Y)=\sqrt{-1} \,\, Y$.
Hence $\g$ has order 4.
\par  b) Each element of $\Gamma$ of order 4  acts on $\bP^1$
with two fixed points  and all other orbits of length 4. So if 
 $\bg$ has  order 4, then  it fixes 2 points in $\P$. Thus $g^2$ has
order 4, by a).
 Hence $g$ has order 8.  

\begin{remark}\label{rem0}
The Lemma implies that an involution of $\bG$ is elliptic if and only if it
fixes no point in its action on the 6-set $P$; equivalently,
if and only if it induces an odd permutation of $P$.
\end{remark}

\begin{remark}\label{rem1}
(i) {\it If a finite subgroup $H$ of $\Gamma$ with $(|H|, char(k))=1$ fixes
a point of $\bP^1$ then $H$ is cyclic:} \ \ 
Indeed, we may assume $H\le \ 
\{\pmatrix 1 & a \\ 0 & b \endpmatrix^\bullet : b\in k^*, a\in k\}$.
The normal subgroup defined by $b=1$ intersects $H$ trivially, 
hence $H$ embeds into its quotient which is isomorphic $ k^*$.
 Hence $H$ is cyclic.

(ii) {\it The degree 2 central extensions of $S_4$:}
\par Their number is $| H^2(S_4,C_2)|=4$ (see
\cite{BS}). 
We construct them as follows.  Let $W$ be the subgroup
 of $GL_4(3)$ generated by
$$
S^\prime=
\pmatrix S & 0 \\
0 & I\\
\endpmatrix, \quad
T^\prime=
\pmatrix T & 0 \\
0 & U \\
\endpmatrix
$$
where $S, T, U\in GL_2(3)= \< S,T \>$ and $ S^3=1=T^2$,
whereas $U$ has order 4.
Then $W$ is a central extension of $PGL_2(3)\iso S_4$ with kernel
 $\{1, w_1,w_2, w_3\}$, where
$$
w_1=\pmatrix I & 0\\
0 & -I
\endpmatrix, \quad
w_2=\pmatrix
-I & 0\\
0 & -I
\endpmatrix, \quad
w_3=w_1w_2.
$$
 \par The $W_i=W/\<w_i\>$, $i=1,2,3$ and the split extension comprise
all  degree 2 central extensions of $S_4$. They are inequivalent since
$W_3$  has no elements of order 8 (as opposed to $W_1$ and $W_2$),
whereas transpositions of $S_4$ lift to involutions (resp., elements
of order 4) in $W_1$ (resp., $W_2$). Note that $W_1\iso GL_2(3)$. 
\end{remark}

\begin{remark}\label{rem2}
Suppose $f_1,f_2,f_3$ are quadratic polynomials in $k[z]$ such that
their product has non-zero discriminant. Then there is an involution
in $\Gamma$ switching the two roots of each $f_i$ if and only if
$f_1,f_2,f_3$ are linearly dependent in $k[z]$ (over $k$).
See Cassels \cite{Cassels}, Thm. 14.1.1,  or Jacobi \cite{J}.
\end{remark}

 \begin{lemma}\label{invol} Suppose $\z$ is an elliptic involution of $G$
and $\e$ its image in $\bG$. Let $u,v$ be the parameters associated with the pair
$(K,\e)$ by Theorem \ref{sect1_thm}. \\
(a) There exists an involution $d$ in $G$ such that the group
$H=<d,\z>$ acts transitively on the 6-set $P$ if and only if 

\begin{equation}\label{D6}
4v-u^2+110u-1125\ \ = \ \ 0
\end{equation}

In this case, $<H,\z_0>\iso D_{12}$ acts as $S_3$ (regularly) on $P$.
\\
(b) There exists an involution $d$ in $G$ such that 
$H=<d,\z>$ has an orbit $Q$ of length 4 on $P$ if and only if
 
\begin{equation}\label{D4}
v^2 \ -\ 4u^3 \ \ = \ \ 0
\end{equation}

In this case, $H\iso D_8$ acts as $V_4$  on $Q$.

(c) If neither (a) nor (b) holds then $G\iso V_4$.
\end{lemma}
\proof 
We may assume that $K=K_{s_1, s_2}$ and $\e=\e_{s_1, s_2}$ as in
Lemma \ref{lemma1}. Then $P=\{a,-a,b,-b,c,-c\}$ for $a,b,c\in k$
with \ $abc=1$,\  $a^2+b^2+c^2=s_1$,\ $a^2b^2+a^2c^2+b^2c^2=s_2$.
Plugging this (with $c=\frac{1}{ab}$)
into \eqref{eq2} expresses $u,v$ as rational functions of $a,b$.
Substituting these expressions for $u,v$ in  \eqref{D6}  and \eqref{D4} yields 

\begin{scriptsize}
\begin{equation}\label{eqD6}
\begin{split}
(a^4b^3-a+a^3b+b+6a^2b^2+ab^3-b^4a^3)(a^4b^3+a-a^3b+b+6a^2b^2-ab^3+b^4a^3)\\
 (a^4b^3-a-a^3b+b-6a^2b^2-ab^3-b^4a^3)(a^4b^3+a+a^3b+b-6a^2b^2+ab^3+b^4a^3) &=0 
\end{split}
\end{equation}
\end{scriptsize}
respectively
\begin{small}
\begin{equation}\label{eqD4}
\begin{split}
 (b-1)^2(b+1)^2(b^2+b+1)^2(b^2-b+1)^2(a-1)^2(a+1)^2(a^2+a+1)^2\\
(a^2-a+1)^2(ab-1)^2(ab+1)^2(a^2b^2+ab+1)^2(a^2b^2-ab+1)^2 &=0
\end{split}
\end{equation}
\end{small}

(a) Such $d$ exists (by Lemma \ref{lem0}) 
if and only if there is an involution $\delta\in\Gamma$ 
fixing $P$ but no point in $P$, and no 4-set in $P$ fixed by $\z$.
By Remark \ref{rem2}, the latter is equivalent to the vanishing of
certain determinants expressed in terms of $a,b$. These 
determinants exactly correspond to the factors in \eqref{eqD6}.
This proves the first claim in (a).

Let $\bar H$ the permutation group on the 6-set $P$ induced by $H$.
We know $\bar H$ is  dihedral and transitive, hence is
(regular) $S_3$ or $D_{12}$. But $D_{12}$ is not generated by two
involutions with no fixed points.  This proves (a).

(b) The first claim is proved as in (a), using the factorization
of $ v^2  - 4u^3$ in \eqref{eqD4}. Now $\bar H$ is  dihedral and transitive
on the 4-set $Q$, hence is $V_4$ or $D_8$. But $D_8$ is not generated
by two involutions with no fixed points.  Thus $H\iso V_4$. 
Since $de$ fixes 
 the two points in $P\setminus Q$, it has 
 order 4. The claim follows.

(c) Suppose neither (a) nor (b) holds. Then $\e$
is the only elliptic involution in $\bG$.
 Hence $\e$ is central in $\bG$.
 If $\gamma$ is another involution in $\bG$, it follows
that $\gamma\e$ is elliptic,
contradiction. Thus $\e$ is the only involution in $\bG$.
Hence either $\bG=<\e>$ or $\bG\iso \bZ_6$. The latter case
cannot occur, see the case $m=6$ in the next section.

\subsection{The list of automorphism groups} \label{3.2}

Since $\bG\hookrightarrow S_6$, all elements of $\bG$ have order $\le6$.
For each $m=4,5,6$ with $(p,m)=1$ there is a unique genus 2 field $K$ such that
$\bG$ contains an element of order $m$.
Indeed, we may assume $\gamma:x\mapsto cx$ with $c\in k^*$ of order $m$.
We may further normalize the coordinate $X$ such that $1\in P$.
Then $P$ consists of all powers of $c$ plus $0$
(for $m\le5$) and $\infty$ (for $m=4$).
 Thus $P$ is also invariant under $x\mapsto 1/x$ for $m=4$
and $m=6$. 
For $p=5$ there is also a unique genus 2 field $K$ such that
$\bG$ contains an element of order $5$.

\medskip
 {\bf I. Sporadic cases:} $\bG$ has elements of order $m\ge4$.

\medskip
$m=4$: \ Here $K$ has equation\ $Y^2=X(X^4-1)$, and $\bG \iso S_4$
(resp., $\bG \iso S_5$, acting as PGL${}_2(5)$ on $P\cong\bP^1(\bF_5)$  )
 if $p\ne5$ (resp., $p=5$).
In each case, $\bG$ is transitive on $P$
and has exactly one class of elliptic involutions (corresponding
to the transpositions in $S_4$ resp. $S_5$). The associated value of $(u,v)$ is
$(5^2,-2\cdot 5^3)$. By Remark \ref{rem1} and Lemma \ref{lem0} we have
$$G\ \  \iso\ \ GL_2(3) \ \ \ \mbox{if} \ \ p\ne5$$
and $$G\ \  \iso\ \ 2^+S_5 \ \ \ \mbox{if} \ \ p=5$$
(the degree 2 cover of $S_5$ where transpositions
lift to involutions).

\medskip
$m=6$: \ If $p=5$ then we are back to the previous case because
$S_5$ has an element of order 6.
The case $p=3$ doesn't
occur here. Now assume $p>5$. Then $K$ has equation\ $Y^2=X^6-1$
and $\bG \iso D_{12}$. Thus $\bG$ has two classes of elliptic involutions,
one of them consisting of the central involution.
The two associated values of $(u,v)$
are $(0,0)$ and $(3^2\ 5^2,\ 2\ 3^3\ 5^3)$. (The first corresponds to the central
involution $x\mapsto -x$ of $\bG$).

By Lemma \ref{invol}(b), the inverse image in $G$ of a Klein 4-subgroup
of $\bG$ is $\iso D_8$. It is a Sylow 2-subgroup of $G$. Thus
$$G\ \  \iso\ \  \bZ_3 \sem D_8$$
where elements of order 4 in $D_8$ act on $\bZ_3$ by inversion.

\medskip
$m=5$: Here $p\ne5$ and $K$ has equation\ $Y^2=X(X^5-1)$.
Further,  $\bG\iso\bZ_5$, \ $G\iso\bZ_{10}$.
There are no elliptic involutions in this case.

\medskip
 {\bf II. The 1-dimensional family with $G\cong D_{12}$} 

Here we assume $\bG$ has an element $\gamma$ of order $3$,
 but none of higher order.
Suppose first $p\ne3$.
Then we may assume $\gamma:x\mapsto cx$ with $c\in k^*$ of order $3$;
also  $1\in P$. Then $P=\{1,c,c^2, a,ac,ac^2\}$ for some $a\in k^*$.
The monic polynomials\  
$(z-1)(z-a)$, $(z-c)(z-c^2a)$, $(z-c^2)(z-ca)$ have the same constant
coefficient, hence are linearly dependent. Hence by Remark
\ref{rem0} there is an elliptic involution $\e$ in $\bG$ with
$\e(1)=a$, $\e(c)=c^2a$, $\e(c^2)=ca$. The group $<\e,\gamma>$
is $\cong S_3$, acting regularly on $P$. Hence by Lemma
\ref{invol}(a) the parameters associated  with the pair $(K,\e)$
satisfy \eqref{D6}:
$$4v-u^2+110u-1125\ \ = \ \ 0$$
Intersection of this curve with $\Delta=0$ is the single point
$(9,54)$. Also the
parameter values $(5^2,-2\ 5^3)$ 
and  $(3^2\ 5^2,\ 2\ 3^3\ 5^3)$ from the
previous case are excluded now. (These values
 satisfy \eqref{D6} which is confirmed by the fact that
the corresponding groups $\bG$  contain a regular $S_3$).
 In the present case, $S_3$ is all of
$\bG$, and by Lemma \ref{invol}(a) we have $G\cong D_{12}$.
If $p=3$ then we may assume $\gamma:x\mapsto x+1$,
and $P=\{0,1,2,a,a+1,a+2\}$. As above we see there is an 
elliptic involution $\e$ in $\bG$ with $<\e,\gamma>\cong S_3$.
The rest is as for $p\ne3$ (only that the parameter value $(0,0)$
doesn't occur because it makes $\Delta$ zero).

\medskip
 {\bf III. The 1-dimensional family with $G\cong D_8$} 

In the remaining cases, $\bG$ has only elements of order $\le2$.
Hence $\bG=\{1\}$, $\bZ_2$ or $V_4$. Here we assume 
$\bG\cong V_4$. Then two of its involutions are elliptic. By Lemma \ref{invol}(b)
it follows that $G\cong D_8$ 
and the $u,v$ parameters satisfy 
$$v^2\ \ = \ \ 4u^3$$
Intersection of this curve with $\Delta=0$ consists of the two points
$(9,54)$ and $(1,-2)$. The values $(0,0)$, $(5^2,-2\ 5^3)$ and
$(3^2\ 5^2,\ 2\ 3^3\ 5^3)$ from Case I are excluded. 

\medskip
 {\bf IV. The 2-dimensional family with $G\cong V_4$}

If $\bG\cong \bZ_2$ then its involution $\e$ is elliptic. 
Indeed, we may assume $\e:x\mapsto -x$ and $1\in P$; if $\e$
is not elliptic then $P=\{0,\infty,1,-1,a,-a\}$ and so $\bG$
contains the additional involution $x\mapsto -a/x$. Thus $G\cong V_4$.
By I--III, this case occurs if and only if the pair $(K,\e)$ has 
$u,v$ parameters with
$$(4v-u^2+110u-1125)\ (v^2-4u^3) \ \ \ne \ \ 0$$

\medskip
 {\bf V. The generic case $G\cong \bZ_2$}

This occurs if and only if $K$ has no elliptic involutions and is not
isomorphic to the field \ $Y^2=X(X^5-1)$. The   existence
of elliptic involutions is equivalent to the condition
in Theorem \ref{mainthm_kap3} (in terms of classical invariants).

\medskip
Summarizing:

\begin{theorem}\label{thm_groups} The  automorphism group $G$ of a
genus 2 field $K$ in characteristic $\ne2$ is isomorphic to 
\ $\bZ_2$, $\bZ_{10}$, $V_4$, $D_8$, $D_{12}$, $\bZ_3 \sem D_8$, $ GL_2(3)$, 
or $2^+S_5$. In the first (resp., last) two cases, $G$ has no
(resp., exactly one) class of elliptic involutions; in the other cases,
it has two classes. Correspondingly, $K$ has either 0, 1 or 2
classes (under $G$-action) of degree 2 elliptic subfields; the case
of one class occurs if and only if $K$ has equation\ $Y^2=X(X^4-1)$.
\end{theorem}

It was noted by Geyer \cite{Geyer} and Gaudry/Schost \cite{Gaudry}
that if $G=D_8$ (resp., $D_{12}$) then degree 2 elliptic subfields
in different classes are 2-isogenous (resp., 3-isogenous). 

\section{The locus of genus 2 curves with elliptic involutions}

\subsection{Classical  invariants and the moduli space $\mathcal M_2$}
Consider a binary sextic i.e. homogeneous polynomial
$f(X,Z)$ in $k[X,Z]$ of degree 6:
$$f(X,Z)=a_6 X^6+ a_5 X^5Z+\dots  +a_0 Z^6$$

  {\bf Classical  invariants} of $f(X,Z)$ are the
following homogeneous polynomials in $k[a_0, \dots , a_6]$ 
of degree $2i$, for  $i=1,2,3,5$. 

\begin{scriptsize}
\begin{equation}
\begin{split}
J_2 := & -240a_0a_6+40a_1a_5-16a_2a_4+6a_3^2\\
J_4 := &
48a_0a_4^3+48a_2^3a_6+4a_2^2a_4^2+1620a_0^2a_6^2+36a_1a_3^2a_5-12a_1a_3a_4^2-12a_2^2a_3
a_5+300a_1^2a_4a_6\\
&
+300a_0a_5^2a_2+324a_0a_6a_3^2-504a_0a_4a_2a_6-180a_0a_4a_3a_5-180a_1a_3a_2a_6+4a_1a_4a
_2a_5\\
& -540a_0a_5a_1a_6-80a_1^2a_5^2\\
J_6 := &-a_5^2a_4a_2+1600a_1^3a_5a_4a_6+1600a_1a_5^3a_0a_2-2240a_1^2a_5^2a_0a_6+20664a_0^2a_4a_
6^2a_2\\
&-640a_0a_4a_2^2a_5^2-18600a_0a_4a_1^2a_6^2+76a_1a_3a_2a_4^3-198a_1a_3^3a_2a_6+26a_1a_3a
_2^2a_5^2\\
& +616a_2^3a_5a_1a_6+28a_1a_4^2a_2^2a_5-640a_1^2a_4^2a_2a_6+26a_1^2a_4^2a_3a_5+616a_1a_
4^3a_0a_5\\
&+59940a_0^2a_5a_6^2a_1+330a_0a_5^2a_3^2a_2+8a_2^2a_3^2a_4^2-24a_2^2a_3^2a_5+60a_2^3a_3^
2a_6-119880a_0^3a_6^3            \\
&
+60a_0a_4^3a_3^2-192a_2^3a_0a_6^2-320a_2^4a_4a_6+176a_1^2a_5^2a_3^2+2250a_1^3a_3a_6^2-9
00a_2^2a_1^2a_6^2\\
& -900a_0^2a_5^2a_4^2-10044a_0^2a_6^2a_3^2+162a_0a_6a_3^4-24a_2^3a_4^3-36a_2^4a_5^2
-36a_1^2a_4^4+76a_2^3a\\
&
-320a_1^3a_5^3+484a+492a_0a_4^2a_2a_3a_5+3060a_0^2a_4a_6a_3a_5-468a_0a_4a_3^2a_2a_6
\\
& +3472a_0a_4a_2a_5a_1a_6-876a_0a_4^2a_1a_6a_3+492a_1a_3a_2^2a_4a_6-238a_1a_3^2a_2a_4a_
5\\
& +1818a_1a_3^2a_0a_6a_5 -1860a_1^2a_3a_2a_5a_6-876a_2^2a_0a_6a_3a_5-3a_5-198a_0a_4a_3^3a_5\\
& +330a_1^2a_3^2a_6a_4-18600a_0^2a_5^2a_6a_2+72a_1a_3^4a_5-24a_1a_3^3a_4^2+2250a_0^2a_5^3a_3
\\
& -1860a_1a_4a_0a_5^2a_3+3060a_1a_3a_0a_6^2a_2\\
J_{10} := & a_6^{-1}Res_X ( f, \frac {\partial f} {\partial X} )\\
\end{split}
\end{equation}
\end{scriptsize}
Here  $J_{10}$ is the discriminant of $f$. 
It vanishes if and only if the binary sextic has a multiple linear factor. These 
$J_{2i}$    are invariant under the natural action of $SL_2(k)$ on
sextics. Dividing such an invariant by another one of the same degree
gives an invariant under $GL_2(k)$ action. 

\par Two genus  2 fields $K$ (resp., curves) in the standard form $Y^2=f(X,1)$ are
isomorphic if and only if the corresponding sextics are $GL_2(k)$
conjugate. Thus if $I$ is a $GL_2(k)$ invariant (resp., homogeneous $SL_2(k)$
invariant), then the expression $I(K)$ (resp., the condition $I(K)=0$)
is well defined. Thus the $GL_2(k)$ invariants are functions 
on the  moduli space $\mathcal M_2$ of genus 2 curves. This $\mathcal
M_2$ is an affine variety with coordinate ring 
$$k[\mathcal M_2]=k[a_0, \dots , a_6, J_{10}^{-1}]^{GL_2(k)}={\textit{subring of degree 0 
elements in}}$$
$k[J_2, \dots ,J_{10}, J_{10}^{-1}]$, see Igusa \cite{Ig}.

\subsection{Classical invariants of genus 2 fields with elliptic involutions}

Under the correspondence in Theorem 4 (resp., Remark 5), the
classical invariants of the field $K$ are:   

\begin{equation}
\begin{split} \label{eq_J}
J_2 & =  240+16u\\
J_4 & =  48v +4u^2+1620-504u \\
J_6 & =  -20664u +96v -424u^2+24u^3+160uv+119880\\
J_{10} & =  64(27-18u-u^2+4v)^2
\end{split}
\end{equation}

respectively 

$$J_2= 384- \frac 1 {16}j$$
$$J_4= 2^{-14} j^2 $$
$$J_6 = 2^{-21} j^2 (-3j+53248)$$
$$J_{10}= 2^{-26} j^4$$

{\bf Proof of Remark 1, concluded:}
The latter formulas explicitly define (in homogeneous coordinates)
the map of   $\mathcal M_1 \setminus \{j=0\}$ to   $\mathcal M_2$
from Remark 1. The function $\frac {J_4 J_6} {J_{10}}\in
k[\mathcal M_2]$  (resp., $\frac {J_2 J_4} {J_6}$) 
is a linear function in $j$ if $char (k)\neq 3$ (resp., $char(k)=3$).
Thus the map is an embedding. This completes the remaining part of the
proof of Remark 1.

\begin{theorem} \label{mainthm_kap3}
The locus $\L_2$ of genus 2 fields with elliptic subfields of degree 2 is the
closed subvariety of $\mathcal M_2$ defined by the equation

\begin{scriptsize}
\begin{equation}
\begin{split}\label{eq_L2_J}
8748J_{10}J_2^4J_6^2- 507384000J_{10}^2J_4^2J_2-19245600J_{10}^2J_4J_2^3
-592272J_{10}J_4^4J_2^2 +77436J_{10}J_4^3J_2^4\\
-81J_2^3J_6^4-3499200J_{10}J_2J_6^3+4743360J_{10}J_4^3J_2J_6-870912J_{10}J_4^2J_2^3J_6 
+3090960J_{10}J_4J_2^2J_6^2\\
-78J_2^5J_4^5-125971200000J_{10}^3 +384J_4^6J_6+41472J_{10}J_4^5+159J_4^6J_2^3
-236196J_{10}^2J_2^5-80J_4^7J_2\\
 -47952J_2J_4J_6^4+104976000J_{10}^2J_2^2J_6-1728J_4^5J_2^2J_6+6048J_4^4J_2J_6^2
-9331200J_{10}J_4^2J_6^2 \\
+12J_2^6J_4^3J_6+29376J_2^2J_4^2J_6^3-8910J_2^3J_4^3J_6^2-2099520000J_{10}^2J_4J_6
+31104J_6^5-6912J_4^3J_6^34  \\
-J_2^7J_4^4 -5832J_{10}J_2^5J_4J_6  -54J_2^5J_4^2J_6^2 +108J_2^4J_4J_6^3  
+972J_{10}J_2^6J_4^2+1332J_2^4J_4^4J_6  = & 0
\end{split}
\end{equation}
\end{scriptsize}

The map  $k^2\setminus \{\Delta=0\} \to \L_2$  described in
Theorem \ref{sect1_thm}  is given  (in homogeneous coordinates) by the
formulas \eqref{eq_J}. It  is birational and surjective if
$char(k)\neq 3$. 
\end{theorem}

\proof
The map is surjective by Theorem  \ref{sect1_thm} and its image is
contained in the subvariety of $\mathcal M_2$ defined by
\eqref{eq_L2_J}; the latter is checked simply by substituting the values
of $J_{2i}$ from \eqref{eq_J}. (We found equation \eqref{eq_L2_J}  by
eliminating $u$ and $v$ from equations \eqref{eq_J}; this equation in
different coordinates was also found in \cite{Gaudry}).

\par Conversely assume $K$ is a genus 2 field with equation
$Y^2=f(X)$ whose classical invariants satisfy \eqref{eq_L2_J}. We
have to show that $K$ has an elliptic involution. We may assume  
$$f(X)=X(X-1)(X-a_1)(X-a_2)(X-a_3)$$
by a coordinate change. Expressing the classical invariants of $K$ in
terms of $a_1, a_2, a_3$, substituting this into \eqref{eq_L2_J} and
factoring the resulting equation yields 

\begin{small}
\begin{equation}
\begin{split}\label{L2_factored}
(a_1a_2-a_2-a_3a_2+a_3)^2(a_1a_2-a_1+a_3a_1-a_3a_2)^2(a_1a_2-a_3a_1-a_3a_2+a_3)^2\\
(a_3a_1-a_1-a_3a_2+a_3)^2(a_1a_2+a_1-a_3a_1-a_2)^2(a_1a_2-a_1-a_3a_1+a_3)^2\\
(a_3a_1+a_2-a_3-a_3a_2)^2(-a_1+a_3a_1+a_2-a_3)^2(a_1a_2-a_1-a_2+a_3)^2\\
(a_1a_2-a_1+a_2-a_3a_2)^2(a_1-a_2+a_3a_2-a_3)^2(a_1a_2-a_3a_1-a_2+a_3 a_2)^2 \\
(a_1a_2-a_3)^2 (a_1-a_3a_2)^2  (a_3a_1-a_2)^2  =&\, 0
\end{split}
\end{equation}
\end{small}
$K$ has an elliptic involution if and
only if there is an involution $\e \in PGL_2(k)$ permuting the set
$\{0,1,\infty, a_1, a_2, a_3\}$ fixed point freely. By Remark \ref{rem2},
the latter is equivalent to the vanishing of certain
determinants expressed in terms of $a_1, a_2, a_3$. These determinants
exactly correspond to the factors  in \eqref{eq_L2_J}. This proves that
$\L_2$ is the closed subvariety of $\M_2$ defined by \eqref{eq_L2_J}.

\par It remains to show  the map in the Theorem is birational.
By Theorem \ref{sect1_thm} we know it is bijective on an open subvariety of
$k^2$. This implies that the corresponding function field extension
$k(u,v)/k(\L_2)$ is purely inseparable, hence its degree $d$ is a
power of $p=char (k)$ (or is 1 in characteristic 0). We need to show $d=1$. 
For this we use the functions $$\frac {J_4} {J_2^2}, \quad
\frac {J_2J_4 -3 J_6} {J_2^3},\quad \frac {J_{10}} {J_2^5}$$
in $k(\mathcal M_2)$. The images of these functions in $k(u,v)$ are: 

\begin{equation}
\begin{split}\label{eq_i}
i_1 &= \frac 1 {64} \frac {12v+u^2+405-126u}{(15+u)^2}\\
i_2 &=- \frac 1 {512} \frac {(-1404v+729u^2-3645+4131u-36uv+u^3)}{(15+u)^3}\\
i_3 &= \frac 1 {16384} \frac {(-27+18u+u^2-4v)^2} {((15+u)^5}
\end{split}
\end{equation}

We compute that $u$ satisfies an equation of degree $\leq 3$  over
the field $k(i_1,i_2)$  whose coefficients are not all zero:

\begin{equation}
\begin{split}
(128i_2-48i_1+1)u^3+(5760i_2+117-3312i_1)u^2+(86400i_2\\
-66960i_1-2349)u+432000i_2-421200i_1+10935 & =0
\end{split}
\end{equation}
Thus $d=1$ (since $p > 3$) and this completes the proof.

\begin{remark}
In characteristic 3 one needs to  replace $v$ by $s_1+s_2$ to get a birational parameterization. 
\end{remark}

\section{Action of $Aut(K)$ on degree $n$ elliptic subfields}

In this section we assume char$(k)=0$. Let $k(X)$, $K$, $G$, $\bG$ as in
section \ref{3.1}
and let $p_1,...,p_6$ the 6 places of $k(X)$ ramified in $K$.

\subsection{Elliptic subfields of $K$ of odd degree}

Consider  an elliptic subfield $F$ of $K$ of odd
degree $n=[K:F]\ge7$. We assume the extension $K/F$ is primitive,
i.e., has no proper intermediate fields. The following facts are
well-known (see \cite{Frey}, \cite{FKV}): The hyperelliptic involution
of $K$ fixes $F$, hence $[F:k(Z)]=2$, where $k(Z)= F\cap k(X)$.
Let $q_1,...,q_r$ be the places of $k(Z)$ ramified in $k(X)$.
Then $r=4$ or $r=5$, and we can label $p_1,...,p_6$  such that
the following holds: $p_i$ lies over $q_i$ for $i=1,2,3$,
and $p_4,p_5,p_6$ lie over $q_4$. Further one of the following holds:

\begin{description}
\item[(1):]  Here $r=5$.
All places of $k(X)$ over $q_1,...,q_4$ different from
$p_1,...,p_6$  have ramification index 2; the $p_i$'s 
have index 1. 
Finally, there is one place $p^{(2)}$ of ramification index 2 over
$q_5$, and all 
other places over $q_5$ have index 1.

\item[(2):] Here and in the following cases we have $r=4$.
Here there is  one place $p^{(4)}$ of ramification index 4 over
$q_4$. All other places of $k(X)$ over $q_1,...,q_4$ different from
$p_1,...,p_6$  have ramification index 2; the $p_i$'s 
have index 1. 

\item[(3):]  Like case (2), only that $p^{(4)}$ lies over $q_1$.

\item[(4):] All places of $k(X)$ over $q_1,...,q_4$ different from
$p_1,...,p_6$  have ramification index 2. The $p_i$'s 
have index 1 except for $p_1$, which has index 3.

\item[(5):]  Like case (4), only now $p_4$ has index 3.

\end{description}

\subsection{Elliptic subfields of $K$ fixed by an automorphism of $K$}
\label{necess}

{\sl Let $g\ne1$ in $\bG=\AA(K)$.
Suppose $g$ fixes $F$. (This is a well-defined statement because
the hyperelliptic involution --- generating the kernel of $G\to\bG$ ---
fixes $F$). Then $g$ has order 2 or 3. 
If $g$  has order 2 it is not an
elliptic involution, and either we are in case (4) and $n\equiv 3$ mod 4, or
we are in case (5) and  $n\equiv 1$ mod 4.
If $g$  has order 3  then either we are in case (1) and $n\not\equiv 1$ mod 3, or
we are in case (2) and  $n\not\equiv 2$ mod 3.} 
 
\medskip\noindent
{\bf Proof:} $g$ acts on $k(X)$ and $k(Z)$,
permuting the ramified places of the extension $k(X)/k(Z)$.
 Thus $g$ fixes the sets
$\{p_1,p_2,p_3\}$ and $\{p_4,p_5,p_6\}$, and the places $p^{(2)}$ resp. $p^{(4)}$.
 Thus $g$ cannot have order
$>3$. Suppose $g$ has order 2. Then it fixes two of the $p_i$'s, hence is
not an elliptic involution and there is no $p^{(2)}$ or $p^{(4)}$.
Thus we are in case (4) or (5). In case (4) (resp., (5)),
 $g$ permutes the $(n-3)/2$ (resp., $(n-5)/2$)
places over $q_1$ (resp., $q_4$) of index 2 fixed point freely,
 hence $n\equiv 3$ mod 4 (resp., $n\equiv 1$ mod 4).
 
Now suppose $g$  has order 3. Then $g$ permutes $p_1,p_2,p_3$
(resp., $p_4,p_5,p_6$) transitively, hence we are in case (1) or (2).
In case (1) (resp., (2)), $g$ fixes $p^{(2)}$ (resp., $p^{(4)}$), 
hence permutes the $n-2$ (resp., $(n-7)/2$)  places
over $q_5$ (resp., $q_4$) of index 1 (resp., 2);
 since it fixes at most one of those places,
we have $n\not\equiv 1$ mod 3 (resp., $n\not\equiv 2$ mod 3).

\subsection{Application of Riemann's existence theorem}

\def\z{\zeta}  \def\s{\sigma}
Let $\z_3$ be a primitive third root of 1 in $k$.
Let $g$ and $F$ as above. 
We can choose the coordinate $Z$ such that
$g(Z)=\z Z$, where $\z=\z_3$ (resp., $\z=-1$) in cases (1) and (2)
(resp., (4) and (5)).  We can further normalize $Z$
such that in case (1) (resp., (2) resp., (4) resp., (5)) 
the places $q_1,...,q_r$ have $Z$-coordinates
$\z^2,1,\z,0,\infty$ (resp., $\infty,1, \z, \z^2$ resp., $0,\infty, 1,-1$
resp., $0,\infty, 1,-1$).

As used in \cite{FKV},
by Riemann's existence theorem 
the equivalence classes of primitive extensions $k(X)/k(Z)$ of degree $n$
with fixed branch points $q_1,...,q_r$ 
and ramification behavior as in (1)--(5) 
correspond to classes of tuples $(\s_1,...,\s_r)$ generating the symmetric group
$S_n$ or alternating group $A_n$ such that $\s_1\cdots\s_r=1$ and

\begin{description}
\item[(1):] \  $\s_i$ is an involution with exactly one fixed point
for $i=1,2,3$, resp., three fixed points for $i=4$, and $\s_5$ is a
transposition.

\item[(2):]\  $\s_i$ is an involution with exactly one fixed point
for $i=1,2,3$, and $\s_4$ has three fixed points, one 4-cycle and the rest
are 2-cycles.

\item[(3):]\  $\s_i$ is an involution with exactly one fixed point
for $i=2,3$, and with three fixed points for $i=4$; and
$\s_1$ has one fixed point, one 4-cycle and the rest
are 2-cycles.

\item[(4):]\ $\s_i$ is an involution with exactly one fixed point
for $i=2,3$, and with three fixed points for $i=4$; and
$\s_1$ has no fixed points, one 3-cycle and the rest
are 2-cycles.

\item[(5):]\   $\s_i$ is an involution with exactly one fixed point 
for $i=1,2,3$, and $\s_4$ has two fixed points, one 3-cycle and the rest
are 2-cycles.

\end{description}

By "classes of
tuples" we mean orbits under the action of $S_n$ by inner automorphisms
(applied component-wise to tuples).
In the case  $k=\bC$, the above correspondence depends on the choice of a
"base point" $q_0$ in $\bP^1\setminus\{q_1,...,q_r\}$ and standard generators 
$\gamma_1,...,\gamma_r$
of the fundamental group $\Gamma(q_0):=\pi_1(
\bP^1\setminus\{q_1,...,q_r\},q_0)$. In particular, $\gamma_1\cdots\gamma_r=1$.
As "base point" we can take any simply connected subset of 
$\bP^1\setminus\{q_1,...,q_r\}$.
The corresponding extensions $\bC(X)/\bC(Z)$ are defined over $\bar\bQ$,
and so one can immediately  pass to the case of general $k$ 
(algebraically closed of char. 0). Here is our choice of the $\gamma_i$
in case (1); we
depict them together with their images $\gamma_i'$ under the map
$z\mapsto\z z$. We depict $\gamma_1,...,\gamma_4$, then $\gamma_5$
is given by the basic relation $\gamma_1\cdots\gamma_5=1$.
All loops are oriented counter-clockwise.

\begin{figure}[ht!]
$$
\entrymodifiers={++[o]}
\xymatrix {
&*++[o][F-]{\zeta \bullet }   &        &  &   
  {}\save[]+<4cm,-1cm>*\txt<8pc>{%
$\gamma_1^\prime=\gamma_2$\\
\, \,     \\
$\gamma_2^\prime=\gamma_3$\\
\, \, \\
$\quad \quad \,\,     \gamma_3^\prime=\gamma_4 \gamma_1 \gamma_4^{-1}$\\
\, \, \\
$\gamma_4^\prime =\gamma_4$\\
\, \, \\
$\quad \quad \,\,     \gamma_5^\prime=\gamma_1^{-1} \gamma_5 \gamma_1$\\
} \restore & 
              \\
& & & &*+[F]{\zeta q_0 \, \bullet} \ar@/_1pc/[ulll]|-{\gamma_2^\prime} \ar@/^/[dr]|-{\gamma_1^
\prime}  
\ar@/^/[dl]|-{\gamma_4^\prime}
\ar@/_2pc/[dddlll]|-{\gamma_3^\prime}  &  \\
& & &*++[o][F-]{0 \bullet}  & &*++[o][F-]{1 \bullet}     \\ 
& & & &*+[F]{q_0 \bullet} \ar@/^1pc/[dlll]|-{\gamma_1} \ar@/_/[ur]|-{\gamma_2} \ar@/^/[ul]|-{
\gamma_4}  
\ar@/_2pc/[uuulll]|-{\gamma_3}   \ar@{~}[uu]|-{Q_0}
    &     \\
&*++[o][F-]{\zeta^2\bullet }   &  & &  \\
}
$$
\caption{The case $q_1,...,q_r\ = \ \z^2,1, \z, 0,\infty$, where $\z=\z_3$ }
\label{fig}
\end{figure}
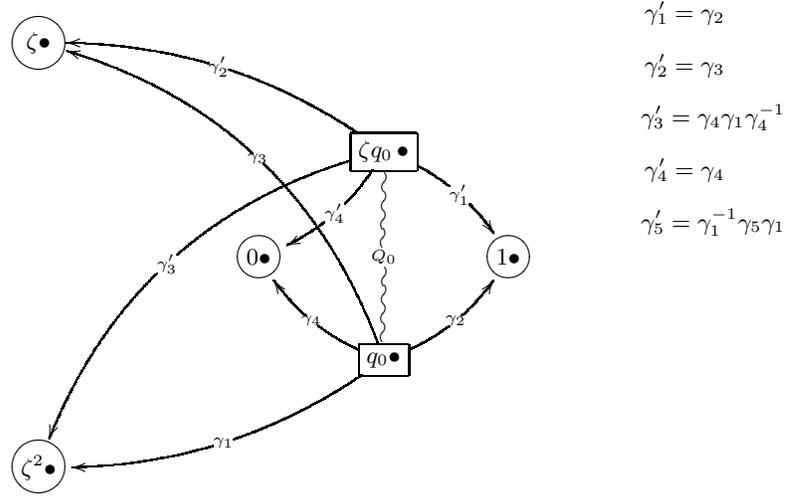

Here we choose $q_0$ as depicted. Let $Q_0$ be the line segment joining
$q_0$ and $\z q_0$. We identify $\Gamma(q_0)$ and $\Gamma(\z q_0)$ via
 the canonical isomorphisms $\Gamma(q_0)\ \cong \  \Gamma(Q_0)\ \cong \ 
\Gamma(\z q_0)$. This yields the above formulas expressing the $\gamma_i'$
in terms of the $\gamma_i$. 

The tuples $(\s_1,...,\s_r)$ corresponding to the extension $\bC(X)/\bC(Z)$,
where $Z=\phi(X)$, are now obtained as follows (see e.g., \cite{Buch},
Ch. 4): Let $\phi$ also denote the map
$\bP^1\to\bP^1$, $x\mapsto \phi(x)$. 
Then lifting of paths gives an action of 
$\Gamma(q_0)$  on $\phi^{-1}(q_0)$, 
hence a homomorphism
of  $\Gamma(q_0)$ to $S_n$. (This homomorphism is determined up to composition
by an inner automorphism of $S_n$ --- re-labeling of the $n$ elements
of $\phi^{-1}(q_0)$ ).  Finally, take $\s_i$ to be the image of $\gamma_i$
under this homomorphism. 

This correspondence between tuples and extensions of $\bC(Z)$ depends 
also on the choice of the coordinate $Z$ (but not on the choice of $X$).
If we replace $Z$ by $Z':=\z Z$,
then the tuple $(\s_1,...,\s_r)$ gets replaced by $(\s_1',...,\s_r')$, where
$\s_i'$ is given in terms of $\s_1,...,\s_r$ by the same formula that
expresses $\gamma_i'$ in terms of $\gamma_1,...,\gamma_r$; see Figure 1
above in case (1). In the other cases (where $r=4$) these formulas appear
already in \cite{Ma} and \cite{Malle}.

\begin{description}
\item[(1)]  { 
\begin{equation}   \s_1^\prime=\s_2\label{formulas} \end{equation}
$$\s_2^\prime=\s_3$$
$$\quad \quad \,\,     \s_3^\prime=\s_4 \s_1 \s_4^{-1}$$
$$\s_4^\prime =\s_4$$
$$\quad \quad \,\,    \s_5^\prime=\s_1^{-1} \s_5 \s_1$$
}
\item[(2)]
 $$\s_1^\prime=\s_2 $$
$$\s_2^\prime=\s_3$$
$$\s_3^\prime=\s_1$$
$$\quad \quad \,\,      \s_4^\prime =\s_1^{-1} \s_4 \s_1$$

\item[(4) and (5)] 
$$\quad \quad \, \,   \s_1^\prime=\s_2 \s_3 \s_2^{-1}$$
$$\s_2^\prime=\s_2$$
$$\s_3^\prime=\s_1$$
$$\quad \quad \, \,  \s_4^\prime =\s_1^{-1} \s_4 \s_1$$
\end{description}

\medskip

Since $Z'=g(Z)=g(\phi(X))=\phi(g(X))$, where $g(X)$ is another
generator of $\bC(X)$, we see that the tuple 
$(\s'_1,...,\s'_r)$  is in the same class as $(\s_1,...,\s_r)$.
Conversely, 
the latter condition is also sufficient for the automorphism $Z\mapsto
\z Z$ to extend to an automorphism of $\bC(X)$. It will permute $p_1,...,p_6$,
hence extend to an automorphism of $K$ fixing $F$.

\subsection{Symmetric tuples}

Primitive extensions $K/F$, where $K$ is a genus 2 field and
and $F$ an elliptic subfield of odd degree $n\ge7$
with fixed branch points of $k(X)/k(Z)$ 
 correspond to classes
of tuples $(\s_1,...,\s_r)$ generating $S_n$ or $A_n$ with $\s_1\cdots\s_r=1$
 as in (1)---(5). Let ${\mathcal T}_j(n)$ be the set of such tuple classes
in case (j), $j=1,...,5$. The number of  these tuple classes grows
polynomially with $n$. (Kani has an exact formula, proved through
a different interpretation of this number, see \cite{K1}). E.g., for $n=7,9,11, 13$ 
we have $|{\mathcal T}_1(n)|= 168, 432, 1100$ and $2184$, respectively.

The condition that $F$ is fixed by an
automorphism of $K$ (different from the identity and
the hyperelliptic involution)
means that $(\s_1,...,\s_r)$ is in the same class as the tuple $(\s'_1,...,\s'_r)$
defined in \eqref{formulas}. Call such tuples {\bf symmetric}. 
 Let ${\mathcal S}_j(n)$ be the set of symmetric tuple classes in
${\mathcal T}_j(n)$. 
The set ${\mathcal S}_j(n)$ can be parameterized  by certain triples,
which we describe in the next section. This allows us to compute the cardinality of
${\mathcal S}_j(n)$ for $n\leq 21$, using a random search to find the triples and 
the structure constant formula \cite{MM}, Prop. 5.5.  to show that we have found all. 
This is based on GAP \cite{GAP} and in particular  \cite{KSH}. Here is the result.

\begin{table}[!hbt]
\renewcommand\arraystretch{1.5}
\noindent\[
\begin{array}{|c|c|c|c|c|c|c|c|c|}
\hline
 & n=7 & n=9 & n=11& n=13 & n=15 & n=17 & n=19& n=21\\
\hline
j=1  & - & 3  & 2  & -  & 6  & 3  & - &2 \\
\hline
j=2 & 1  &  0 & - & 2  & 0  & - & 4 & 0 \\
\hline
j=4 & 2  & - & 3  & - & 4 & - & 5 & -\\
\hline
j=5 & - & 3  &- & 3  &- & 4 & - & 5\\
\hline
\end{array}
\]
\caption{ $|\mathcal S_j(n) |$ = number of symmetric tuple classes}
\label{ }
\end{table}

From the table it appears that the necessary conditions in section \ref{necess}
(for the existence of extensions $K/F$ with non-trivial automorphisms)
are sufficient in most cases (at least for those $n$ in reach of computer
calculation). It is intriguing that the number of these extensions seems
to be very small, but mostly $>1$.

\subsection{Parameterization of symmetric tuples}

\def\t{\tau} \def\r{\rho}  \def\I{\mbox{Ind}}

Let $(\s_1,...,\s_5)$ be a tuple representing an element of  ${\mathcal S}_1(n)$.
Thus there is $\t\in S_n$ with $\s'_i\ = \ \s_i^\t$ for $i=1,...,5$.
Then $\s_i^{\t^3}\ = \ \s_i^{\s_4}$, hence $\t^3=\s_4$.
Thus all $\s_i$ can be expressed in terms of $\t$ and $\s:=\s_1$:
\begin{equation}
\s_1\ = \  \s, \ \ \ \s_2\ = \  \s^\t,\ \ \ \s_3\ = \  \s^{\t^2},
\ \ \ \s_4\ = \ \t^3,\ \ \ \s_5\ = \  (\s\t^{-1})^3
\label{5br_pt}
\end{equation}

Passing from $(\s,\t\, \r)$
to $(\s_1,...,\s_5)$ is a case of "translation", see \cite{Magaard}
and \cite{Malle}.
Recall that the {\bf index} $\I(\pi)$ of $\pi\in S_n$ is defined as
$n$ minus the number of orbits of $\pi$. Since $\s=\s_1$ is an
involution with exactly one fixed point, we have  $\I(\s)=(n-1)/2$.
From $\t^3=\s_4$ it follows that

\begin{equation}
\I\ (\r)\ \ \le 
\left\{ 
\aligned
\frac{5(n-3)}{6}\  +\ 2 \ \ \ \  \ \mbox{if}\ \ n\equiv0 \mod 3 \\
\frac{5(n-5)}{6}\  +\ 3 \ \ \ \  \ \mbox{if}\ \ n\equiv2 \mod 3 
\endaligned
\right.
\end{equation}
where equality holds if and only if $\t$ has cycle type as in the Lemma below (case $j=1$).
Further, for $\r:=\s\t^{-1}$ we have $\rho^3=\s_5$ (a transposition).
Hence
\begin{equation}
\I\ (\r)\ \ \le 
\left\{ 
\aligned
\frac{2(n-3)}{3}\  +\ 1 \ \ \ \  \ \mbox{if}\ \ n\equiv0 \mod 3\\
\frac{2(n-2)}{3}\  +\ 1 \ \ \ \  \ \mbox{if}\ \ n\equiv2 \mod 3
\endaligned
\right.
\end{equation}
where equality holds if and only if $\r$ is as in the Lemma below (case $j=1$).

It follows that $\I(\s)+\I(\t)+\I(\r)\ \le\ \ 2(n-1)$. 
The reverse inequality holds by the Riemann Hurwitz formula since
$<\s,\t,\r>=S_n$. Hence $\t$ and $\r$ are
of cycle type as claimed in the following Lemma.

\begin{lemma}\label{symm2} There is a bijection between
 ${\mathcal S}_j(n)$ and the set of classes of triples $(\s,\t,\r)$
generating $S_n$ (resp., $A_n$) with  $\r\t=\s$, where
$\s$ is an involution with exactly one fixed point and $\t,\r$ are
of the following cycle type:\\
\begin{description}
\item[j=1:] \  $\r$ has one 2-cycle, at most one fixed point
and the rest are 3-cycles;\\ 
\quad  \quad  $\t$ has one 3-cycle, at most one 2-cycle and the rest are
6-cycles.

\item[j=2:]\ $\t$ has at most one fixed point and its other cycles are
all 3-cycles; \\
\quad \,  \quad  $\r$ has one 4-cycle, one 3-cycle, at most one 2-cycle and the rest are
6-cycles.

\item[j=4:]\ $\r$ has one fixed point, one 2-cycle and the rest are
4-cycles;  \\
\quad \,        $  \t$ has one 3-cycle and the rest are
4-cycles.

\item[j=5:]\  $\r$ has one 2-cycle, one 3-cycle and the rest are
4-cycles;  \\
\quad \,  \quad        $\t$ has one fixed point and its other cycles are
all 4-cycles. 

\end{description}
\end{lemma}
\proof We only discuss case (1), the other cases are similar.
In this case, it remains to show that for  given $\s,\t,\r$ as in the Lemma,
 formulas \eqref{5br_pt}
define a tuple $(\s_1,...,\s_5)$ representing an element of ${\mathcal S}_1(n)$.
First one verifies that the tuple $(\s'_1,...,\s'_5)$ defined as in
\eqref{formulas} is conjugate to $(\s_1,...,\s_5)$ under $\t$.
This implies  that $<\s_1,...,\s_5>$ is normal in
$<\s,\t>=S_n$, hence equals $S_n$ (since it contains a transposition).

\end{document}